\documentclass[11pt]{article}
\usepackage{epic,latexsym,amssymb,url}
\usepackage{color}
\usepackage{tikz}

\usepackage{geometry,graphicx}

\newtheorem{theorem}{Theorem}

\newtheorem{claim}{Claim}
\newtheorem{corollary}{Corollary}
\newtheorem{observation}{Observation}

\newtheorem{conjecture}{Conjecture}

\newcommand{\QED}{$\Box$}
\newcommand{\smallqed}{{\tiny ($\Box$)}}

\newcommand{\diam}{{\rm diam}}

\newcommand{\D}{{Dominator }}
\newcommand{\St}{{Staller }}

\newcommand{\igt}{\iota_{\rm gt}}
\newcommand{\igtS}{\iota_{\rm gt}'}

\newenvironment{unnumbered}[1]{\trivlist \item [\hskip \labelsep {\bf
#1}]\ignorespaces\it}{\endtrivlist}

\newcommand{\proof}{\noindent\textbf{Proof. }}
\newcommand{\2}{ \vspace{0.2cm} }
\newcommand{\1}{ \vspace{0.1cm} }

\let\oldenumerate\enumerate
\renewcommand{\enumerate}{
  \oldenumerate
  \setlength{\itemsep}{0pt}
  \setlength{\parskip}{0pt}
  \setlength{\parsep}{0pt}
}

\begin{document}

\title{Total isolation game in graphs}

\author{$^1$Michael A. Henning and $^2$Douglas F. Rall
\\ \\
$^1$Department of Mathematics and Applied Mathematics \\
University of Johannesburg \\
Auckland Park, 2006 South Africa\\
\small \tt Email: mahenning@uj.ac.za  \\
\\
$^2$Professor Emeritus of Mathematics \\
Furman University \\
Greenville, SC, USA\\
\small \tt Email: doug.rall@furman.edu}

\date{}
\maketitle

\begin{abstract}
The total isolation game is played on a graph $G$ by two players  who take turns playing a vertex such that if $S$ is the set of already played vertices, then a vertex can be selected only if it is adjacent to a vertex that belongs to a (nontrivial) component of the graph $G - N_G(S)$ of order at least~$2$ or a vertex that is isolated in $G - N_G(S)$ and belongs to the set $S$, where $N_G(S)$ is the set of vertices adjacent to a vertex in $S$. Dominator wishes to finish the game with the minimum number of played vertices, while Staller has the opposite goal. The game total isolation number $\igt(G)$ is the number of moves in the Dominator-start game where both players play optimally. We prove that if $G$ is a connected graph of order~$n \ge 3$, then $\igt(G) < \frac{5}{6}n$. Furthermore if $G$ has minimum degree at least~$2$, then we prove that $\igt(G) \le \frac{3}{4}n$. More generally, if $G$ is a connected graph of order~$n \ge 3$ with minimum degree~$\delta$ where $\delta \ge 2$, then we prove that $\igt(G) \le \left( \frac{2\delta-1}{3\delta-2} \right) n$. Among other results it is proved that if $G$ is a graph of order~$n$ with diameter~$2$, then $\igt(G) \le \frac{2}{3}n$.
\end{abstract}


\section{Introduction}

A vertex $u$ in a graph $G$ \emph{dominates} a vertex~$v$ if $u = v$ or $u$ is adjacent to $v$ in $G$. A \emph{dominating set} of $G$ is a set $S$ of vertices of $G$ such that every vertex in $G$ is dominated by a vertex in $S$. A vertex $u$ in a graph $G$ \emph{totally dominates} a vertex~$v$ if $u$ is adjacent to $v$ in $G$. A \emph{total dominating set} of $G$ is a set $S$ of vertices of $G$ such that every vertex of $G$ is totally dominated by a vertex in $S$. Thus, a set $S$ is a dominating set of $G$ if every vertex in $V(G) \setminus S$ is adjacent to at least one vertex in $S$, and a set $S$ is a total dominating set of $G$ if every vertex is adjacent to at least one vertex in $S$. The concept of domination and its variations have been widely studied in theoretical, algorithmic and application aspects. A thorough treatment of (total) domination in graphs and its variants can be found in the books~\cite{HaHeHe-23,HeYe-book}.

The \emph{open neighborhood} of a vertex $v$ in $G$ is the set of neighbors of $v$, denoted $N_G(v)$, where two vertices are \emph{neighbors} if they are adjacent. Thus, $N_G(v) = \{u \in V(G) \, \colon uv \in E(G)\}$. The \emph{closed neighborhood} of $v$ is the set $N_G[v] = N_G(v) \cup \{v\}$. The \emph{open neighborhood} of a set $S \subseteq V(G)$ is the set of all neighbors of vertices in $S$, denoted $N_G(S)$, whereas the \emph{closed neighborhood} of $S$ is $N_G[S] = N_G(S) \cup S$. Thus,
\[
N_G(S) = \bigcup_{v \in S} N_G(v) \hspace*{0.5cm} \mbox{and} \hspace*{0.5cm} N_G[S] = \bigcup_{v \in S} N_G[v].
\]

A set $S \subseteq V(G)$ is an \emph{isolating set} of $G$ if removing $S$ and its neighborhood $N_G(S)$ leaves no edge, that is, $V(G) \setminus N_G[S]$ is an independent set in $G$. The concept was introduced (in a more general setting) in 2017 by Caro and Hansberg~\cite{CaHa-17}. The \emph{isolation number} $\iota(G)$ of $G$  is the minimum cardinality among all isolating sets of $G$. This parameter is equivalent to the vertex-edge domination number, introduced earlier by~Peters~\cite{Pe-86}. Isolation in graphs is studied further, for example, in~\cite{Bo-25,BoGo-24,CaAnWu-25,DaLeSoVa-21,GoHe-25,LeMoSo-24,Zy-19}.

A \emph{total isolating set} is an isolating set $S$ with the additional property that the subgraph $G[S]$ induced by $S$ contains no vertex of degree~$0$. The \emph{total isolation number} $\iota^t(G)$ of $G$ is the minimum cardinality among all total isolating sets of $G$. This parameter was originally introduced by Boutrig and Chellali~\cite{BoCh-18} as the \emph{total vertex-edge domination number}, and was subsequently studied, for example, in~\cite{BoGoHe-25a,BoGoHe-25b,CaAnWu-25}.

In 2010 Bre{\v{s}}ar, Klav{\v{z}}ar, and Rall~\cite{BrKlRa-10} published the seminal paper on the domination game which belongs to the growing family of competitive optimization graph games. Domination games played on graphs are now very well studied in the literature. The subsequent rapid growth by the scientific community of research on domination games played on graphs inspired the recent book entitled ``Domination games played on graphs" by Bre\v{s}ar, Henning, Klav\v zar, and Rall~\cite{BrHeKlRa-21}, which presented the state of the art results at the time and shows that the area is rich for further research. In this paper, we study the total version of the isolation game in graphs.

Recently Bre\v sar, Dravec, Johnston, Kuenzel, and Rall in~\cite{BrDrJoKuRa-24}, introduced a game counterpart of the isolation number, called the  \emph{isolation game}, following the ideas of the classical domination game~\cite{BrKlRa-10}. The isolation game was subsequently studied further in~\cite{BuDrHeKl-17} at the Workshop on Games on Graphs III, held in Rogla, Slovenia, in June 2025. In this paper, we introduce and study a game counterpart of the total isolation number.

The \emph{total isolation game} is played on a graph $G$ by Dominator and Staller, who take turns selecting/playing a vertex from $G$ while obeying the rule that if $S$ is the set of already played vertices, then a vertex~$v$ can be selected only if $v$ totally dominates a vertex~$u$ that belongs to a (nontrivial) component of the graph $G - N_G(S)$ of order at least~$2$ or $v$ totally dominates a vertex $u\in S$ that is isolated in $G - N_G(S)$.  Such a vertex $v$ is \emph{playable}, and a move selecting~$v$ is \emph{legal}. The game ends when no playable vertex exists. When the game ends, the set of vertices selected forms a total isolating set of $G$.

As in the domination game, in the total isolation game Dominator wishes to finish the game by playing the fewest number of vertices, while Staller wishes to delay the process as much as possible. Thus, Dominator seeks to minimize the size of the chosen set while Staller tries to make it as large as possible. If Dominator starts the game, we speak of a \emph{D}-\emph{game}, otherwise it is an \emph{S}-\emph{game}. The \emph{game total isolation number}, $\igt(G)$, and the \emph{Staller}-\emph{start game total isolation number}, $\igtS(G)$, are the number of vertices selected in the D-game and the S-game, respectively, provided that both players play optimally.

For graph theory notation and terminology, we generally follow~\cite{HaHeHe-23}.  Specifically, let $G$ be a graph with vertex set $V(G)$ and edge set $E(G)$, and of order $n(G) = |V(G)|$ and size $m(G) = |E(G)|$. We denote the degree of a vertex $v$ in $G$ by $\deg_G(v)$. The minimum and maximum degrees in $G$ are denoted by $\delta(G)$ and $\Delta(G)$, respectively. An \emph{isolated vertex} in $G$ is a vertex of degree~$0$ in $G$. A graph is \emph{isolate}-\emph{free} if it contains no isolated vertex. A set $X$ is a \emph{packing} in $G$ if for any two distinct vertices $u,v \in X$, we have $d_G(u,v) \ge 3$, and so $N_G[u] \cap N_G[v] = \emptyset$.

A \emph{trivial graph} is the graph of order~$1$, and a \emph{nontrivial graph} has order at least~$2$. A \emph{cycle} on~$n$ vertices is denoted by $C_n$ and a \emph{path} on $n$ vertices by $P_n$. The  \emph{complete graph} on $n$ vertices is denoted by $K_n$. For a subset $S$ of vertices of a graph $G$, we denote by $G - S$ the graph obtained from $G$ by deleting the vertices in $S$ and all edges incident with vertices in $S$. If $S = \{v\}$, then we simply write $G - v$ rather than $G - \{v\}$. The subgraph induced by the set $S$ is denoted by $G[S]$. The diameter of $G$, denoted $\diam(G)$, is the maximum distance among all pairs of vertices in~$G$.

\section{Continuation Principle}

One of the main tools in analyzing the domination game is the Continuation Principle first presented by Kinnersley, West, and Zamani~\cite{KiWeZa-13}. An important consequence of the Continuation Principle for a given game variant of the domination game is the fundamental property that the number of moves in the D-game and the S-game when played optimally can differ by at most~$1$. It is known, for example, that the Continuation Principle holds for the total domination game. And as a special case of a more general result, the authors in~\cite{BrDrJoKuRa-24} proved that the Continuation Principle holds in the isolation game.

Consider the total isolation game played on a path $P_5$. Let $G$ be the path $P_5$ given by $v_1v_2v_3v_4v_5$. In the D-game, \D plays as his first move the central vertex~$v_3$, and so the set of played vertices immediately after Dominator's first move is $S = \{v_3\}$. We note that no component of $G - N_G(S)$ is nontrivial and that the vertex $v_3$ is the only isolated vertex in $G - N_G(S)$ that belongs to the set $S$. Hence the only legal move available to \St is either to play the vertex~$v_2$ or the vertex~$v_4$ in order to totally dominate the vertex~$v_3$. Once \St plays her first move (either $v_2$ or $v_4$), the game is complete. Hence, $\igt(G) = 2$. In the S-game, \St plays as her first move the leaf~$v_1$. If \D plays~$v_2$ on his first move, then \St responds by playing~$v_5$. If \D plays~$v_3$ on his first move, then \St responds by playing~$v_4$. If \D plays~$v_4$ on his first move, then \St responds by playing~$v_5$. If \D plays~$v_5$ on his first move, then \St responds by playing~$v_4$. In all four cases, one additional move is required by \D to complete the game. Hence, $\igtS(G) \ge 4$. (One can show that four moves always suffice to complete the S-game, implying that $\igtS(G) = 4$, but we do not need this stronger property here.) Therefore, $|\igt(G) - \igtS(G)| \ge 2$.

Since the total isolation game combines the flavor of the total domination game and the isolation game, and since the Continuation Principle holds for both these games, one might expect that the Continuation Principle holds for the total isolation game. If this were the case, then for every isolate-free graph $G$ we would have the desirable property that $|\igt(G) - \igtS(G)| \le 1$. However, since there are graphs $G$ for which $|\igt(G) - \igtS(G)| > 1$, we infer the following result.

\begin{observation}
\label{Cont-Principle}
The Continuation Principle does not hold for the total isolation game.
\end{observation}

\section{Upper bounds for given minimum degree}
\label{S:min-degree-large}

We remark that the Continuation Principle is a powerful tool for establishing upper bounds on game domination type parameters, including game domination number, the game total domination number, the game isolation number, to name games where the Continuation Principle holds. However for several variants of the domination games, such as the independent domination game and the connected domination game, in which the Continuation Principle does not hold, it has proved challenging to find good upper bounds on the associated game parameters. Observation~\ref{Cont-Principle} therefore suggests that finding good upper bounds on the game total isolation number is likely to be a difficult problem. In this section, our aim is to provide upper bounds on the game total isolation number in terms of the minimum degree of the graph.

During the course of the total isolation game, if $S$ is the set of selected vertices played to date, then a vertex $v$ is considered \emph{marked} if $v \in N_G(S)$ or if $v \notin S$ and $v$ is isolated in $G - N_G(S)$. We note that a marked vertex may or may not be playable. A vertex is \emph{unmarked} if it is not marked. At any given stage of the total isolation game, we define the sets $S$, $M$ and $U$ as follows: \\ [-24pt]
\begin{enumerate}
\item[$\bullet$]  $S$ is the set of selected vertices.
\item[$\bullet$]  $M$ is the set of marked vertices.
\item[$\bullet$]  $U$ is the set of unmarked vertices.
\end{enumerate}

We note that if $u$ is an unmarked vertex, then $u \in U$ and either $u \in S$ and $u$ is an isolated vertex in $G - N_G(S)$ or $u$ belongs to a nontrivial component in the graph $G - N_G(S)$. Thus every playable vertex in the total isolation game is adjacent to at least one unmarked vertex, and after such a vertex is played, at least one new vertex is marked.  The game ends when all vertices of $G$ are marked.

\begin{theorem}
\label{thm:delta-bound}
If $G$ is a connected graph of order~$n \ge 3$ with minimum degree~$\delta$ where $\delta \ge 2$ and maximum degree~$\Delta$, then
\[
\igt(G) \le \left( \frac{2\delta-1}{3\delta-2} \right) n - \frac{\Delta - 2}{3\delta-2}.
\]
\end{theorem}
\proof
Let $G$ be a connected graph of order~$n$ with minimum degree~$\delta$ where $\delta \ge 2$ and maximum degree~$\Delta$. As the D-game is played, let $m_{2i-1}$ be the number of new vertices marked by \D on his $i$th move, and let $m_{2i}$ be the number of new vertices marked by \St on her $i$th move. Upon completion of the game, say in $t$ moves, all vertices are marked implying that

\begin{equation}
\label{Eq1}
n = \sum_{i=1}^{\lceil \frac{t}{2} \rceil} m_{2i-1} + \sum_{i=1}^{\lfloor \frac{t}{2} \rfloor} m_{2i}.
\end{equation}

In what follows, let \D play according to a greedy strategy that on each of his moves he plays a vertex that results in a maximum increase in the number of marked vertices. That is, on each of his moves, \D plays a vertex that marks as many previously unmarked vertices as possible. Before Dominator's first move, all vertices are unmarked, and so $U = V(G)$ and $S = M = \emptyset$. Let $t$ denote the number of moves played upon completion of the game when \D adopts his greedy strategy.
We consider two stages in the game.

\paragraph{Stage~$1$.} A move of Dominator belongs to Stage~$1$ if there exists a playable move available to \D that marks at least two new vertices. A move of Staller belongs to Stage~$1$, if so does Dominator's previous move.

\paragraph{Stage~$2$.} A move of Dominator belongs to Stage~$2$ if every playable move available to \D marks exactly one new vertex. A move of Staller belongs to Stage~$2$, if so does Dominator's previous move.

As remarked earlier, \D plays according to a greedy strategy on each of his moves. We consider first the case when the game is complete in Stage~$1$.

\begin{claim}
\label{c:claim-1}
If the game is complete in Stage~$1$, then $t \le \frac{2}{3}n - \frac{2}{3}(\Delta - 2)$.
\end{claim}
\proof
If \D plays a vertex~$v$ of maximum degree~$\Delta$ on his first move, then he marks at least~$\Delta$ vertices, namely all~$\Delta$ vertices in $N_G(v)$ together with all isolated vertices in $G - N_G(v)$ different from~$v$ if such vertices exist. Thus there exists a playable move available to \D that marks at least $\Delta$ vertices. According to Dominator's greedy strategy, we therefore infer that on his first move Dominator marks at least $\Delta$ vertices, and so $m_1 \ge \Delta \ge 2$. Every playable move marks at least one new vertex. In particular, every move of Staller marks at least one new vertex, and so $m_{2i} \ge 1$ for all $i \in \{1, \ldots, \lfloor t/2 \rfloor\}$.

Suppose that Stage~1 consists of $k_1$ moves of Dominator. Dominator's first move marks at least~$\Delta$ vertices and each of his remaining $k_1 - 1$ moves in Stage~1 mark at least two new vertices, and so $m_1 \ge \Delta$ and, if $k_1 \ge 2$, then $m_{2i-1} \ge 2$ for all $i \in [k_1] \setminus \{1\}$. If the game is complete immediately after Dominator plays his $k_1^{\rm th}$ move, then $t = 2k_1 - 1$ and, by Inequality~(\ref{Eq1}) and our earlier observations, we have
\[
n = \sum_{i=1}^{k_1} m_{2i-1} + \sum_{i=1}^{k_1 - 1} m_{2i} \ge (2k_1 + \Delta - 2) + (k_1 - 1) = 3k_1 + \Delta - 3,
\]
and so
\[
t = 2k_1 - 1 = \frac{2}{3}(n - \Delta + 3) - 1 < \frac{2}{3}n - \frac{2}{3}(\Delta - 2), \1
\]
and so the desired upper bound in the statement of the claim holds (with strict inequality). Hence, we may assume that $t \ge 2k_1$. If $t = 2k_1$, then by Inequality~(\ref{Eq1}) we have $n \ge (2k_1 + \Delta - 2) + k_1 = 3k_1 + \Delta - 2$, and in this case
\[
t = 2k_1 \le \frac{2}{3}(n  - \Delta + 2),
\]
and so, once again, the desired upper bound in the statement of the claim holds.~\smallqed

\medskip
Suppose, next, that the game is not complete in Stage~$1$ and that  Stage~$1$ consists of $k_1$ moves of Dominator. Thus, $t > 2k_1$. In Stage~$2$ of the game, every playable move by either player marks exactly one new vertex. Suppose that $k_2$ vertices in total are played in Stage~$2$ of the game, and so
\begin{equation}
\label{Eq3}
t = 2k_1 + k_2.
\end{equation}

By our earlier assumptions, we have $t > 2k_1$, and so $k_2 \ge 1$. We now consider the sets $S$, $M$ and $U$ of selected, marked and unmarked vertices immediately after Stage~1 and before the game enters Stage~2. We proceed further by proving some structural properties that hold in the graph $G$.

\begin{claim}
\label{c:claim-2}
The following properties hold in the graph $G$. \\ [-24pt]
\begin{enumerate}
\item[{\rm (a)}]  $U$ is an independent set in $G$.
\item[{\rm (b)}]  $U \subseteq S$.
\item[{\rm (c)}]  $U$ is a packing in $G$.
\end{enumerate}
\end{claim}
\proof
(a) Suppose, to the contrary, that $U$ is not an independent set, and let $C$ be a nontrivial component in $G[U]$. Suppose that $|V(C)| \ge 3$. In this case, the component $C$ contains a vertex, say~$v$, of degree at least~$2$ in $C$. If \D plays the vertex~$v$, then he marks the neighbors of $v$ in $C$ and possibly some additional vertices not adjacent to~$v$. Such a move of \D therefore marks at least two new vertices, contradicting the fact that the game is in Stage~$2$. Hence, $|V(C)| = 2$, and so $C$ is a $P_2$-component in $G[U]$.

Let $C$ be the path $u_1u_2$. If $u_i \in S$, then the vertex $u_{3-i}$ would be a marked vertex for $i \in [2]$, a contradiction. Hence, neither $u_1$ nor $u_2$ belongs to the set $S$. The graph $G$ is connected and $\Delta \ge 2$. Hence renaming vertices in the component $C$ if necessary, we may assume that $\deg_G(u_1) \ge 2$. Let $u$ be a neighbor of $u_1$ different from~$u_2$. Necessarily, the vertex $u$ is marked but has not been selected.  That is, $u \in M \setminus S$. If \D plays the vertex~$u$, then he marks its neighbor $u_1$ that belongs to~$C$. Furthermore, the vertex~$u_2$ now becomes isolated in the graph $G - N_G(S \cup \{u\})$ if $uu_2 \notin E(G)$, and~$u_2$ becomes marked if $uu_2 \in E(G)$.  In either case, $u_2$ becomes marked if \D plays~$u$. Therefore, \D marks at least two new vertices when he plays the vertex~$u$, contradicting the fact that the game is in Stage~$2$. This proves part~(a).

(b) If there is a vertex, say~$x$, in $U$ that does not belong to the set $S$, then $x \notin S$ and $x$ is isolated in $G - N_G(S)$. However this would imply that the vertex~$x$ is a marked vertex, a contradiction. Therefore, every vertex in $U$ belongs to the set $S$, and so $U \subseteq S$. This proves part~(b).

(c) Suppose, to the contrary, that $U$ is not a packing in $G$. Hence, since $U$ is an independent set, we infer that there exists two distinct vertices $u_1$ and $u_2$ in $U$ that have a common neighbor, say~$v$, in $G$. Necessarily, the vertex~$v$ is marked but has not been selected.  If \D plays the vertex~$v$, then he marks both~$u_1$ and $u_2$ and possibly some additional vertices. Such a move of \D therefore marks at least two new vertices, contradicting the fact that the game is in Stage~$2$. This proves part~(c), and completes the proof of the claim.~\smallqed

\medskip
By Claim~\ref{c:claim-2}(c), the set $U$ is a packing in $G$ and $U \subseteq S$. We infer that every move played in Stage~$2$ marks exactly one vertex that belongs to the set $U$. Let $U = \{u_1,\ldots,u_{k_2}\}$ and let $B_i = N_G(u_i)$ for $i \in [k_2]$. We note that the sets $B_1, \ldots,B_{k_2}$ are pairwise vertex disjoint. We now let $B = N_G(U)$, and so $B = \bigcup_{i=1}^{k_2} B_i$ and
\begin{equation}
\label{Eq4}
|B| = \sum_{i=1}^{k_2} |B_i| = \sum_{i=1}^{k_2} \deg_G(u_i) \ge \delta  k_2.
\end{equation}

If the set $B$ contains a vertex that belongs to the set $S$, then the neighbor of such a vertex that belongs to the set $U$ would be marked, a contradiction. Hence, $B \cap S = \emptyset$.

Recall that by Equation~(\ref{Eq3}), we have $t = 2k_1 + k_2$ where $t$ denote the number of moves played upon completion of the game when \D adopts his greedy strategy. Further recall that we are considering the sets $S$, $M$ and $U$ of selected, marked and unmarked vertices immediately after Stage~1 and before the game enters Stage~2. Thus at this stage of the game, exactly $2k_1$ vertices have been played.

Let $A$ be the set of vertices at distance at least~$2$ from every vertex in~$U$, and so $A = V(G) \setminus N_G[U]$. Recall that $U \subseteq S$ and $|U| = k_2$, and so the set $S$ contains all $k_2$ vertices in $U$. Since $B \cap S = \emptyset$, the $2k_1 - k_2$ vertices in $S$ that do not belong to~$U$ therefore belong to the set~$A$. Hence,
\begin{equation}
\label{Eq5}
|A| \ge |A \cap S| = 2k_1 - k_2.
\end{equation}

Recall that by supposition $\delta \ge 2$. We note that the sets $A$, $B$ and $U$ are vertex disjoint and $V(G) = A \cup B \cup U$. Thus, by Inequalities~(\ref{Eq4}) and~(\ref{Eq5}), and recalling that $|U| = k_2$, we have
\[
\begin{array}{lcl}
k_2 & = & n - |A| - |B| \\ \1
& \le & n - (2k_1 - k_2) - \delta  k_2 \\ \1
& = & n - 2k_1 - (\delta - 1)k_2,
\end{array}
\]
and so $(\delta - 1)k_2 \le n - (2k_1 + k_2) = n - t$, or, equivalently,
\begin{equation}
\label{Eq6}
k_2 \le \frac{n-t}{\delta - 1}.
\end{equation}

Immediately upon completion of Stage~$1$ of the game, we have
\[
|M| = \sum_{i=1}^{k_1} m_{2i-1} + \sum_{i=1}^{k_1} m_{2i} \ge (2k_1 + \Delta - 2) + k_1 = 3k_1 + \Delta - 2.
\]

Recall that $t = 2k_1 + k_2$ by Equation~(\ref{Eq3}), and so $k_1 = \frac{1}{2}(t - k_2)$. From these observations and by Inequality~(\ref{Eq6}), we have
\[
\begin{array}{lcl}
n & = & |M| + |U| \\ \1
& \ge & (3k_1 + \Delta - 2) + k_2 \\ \1
& = & \frac{3}{2}(t - k_2) + k_2 + \Delta - 2 \\ \1
& = & \frac{3}{2}t - \frac{1}{2}k_2 + \Delta - 2 \\ \1
& \ge & \frac{3}{2}t - \frac{n-t}{2(\delta - 1)} + \Delta - 2,
\end{array}
\]
or, equivalently,

\begin{equation}
\label{Eq7}
\left( 1 + \frac{1}{2(\delta - 1)} \right) n \ge \left( \frac{3}{2} + \frac{1}{2(\delta - 1)} \right) t + \Delta - 2.
\end{equation}

Inequality~(\ref{Eq7}) yields the simplified inequality
\begin{equation}
\label{Eq8}
\left( \frac{2\delta - 1}{2\delta - 2} \right) n \ge \left( \frac{3\delta - 2}{2\delta - 2} \right) t + \Delta - 2.
\end{equation}

Solving for $t$ yields the inequality
\begin{equation}
\label{Eq9}
t \le \left( \frac{2\delta-1}{3\delta-2} \right) n - \frac{\Delta - 2}{3\delta-2}.
\end{equation}

Recall that $t$ denotes the number of moves played in the game when \D adopts his greedy strategy. By Claim~\ref{c:claim-1}, if the game is complete in Stage~$1$, then $t \le \frac{2}{3}n - \frac{2}{3}(\Delta - 2)$. If the game is not complete in Stage~$1$ and therefore enters Stage~$2$, then by Inequality~(\ref{Eq9}) we have $t \le \left( \frac{2\delta-1}{3\delta-2} \right) n - \frac{\Delta - 2}{3\delta-2}$. For $n \ge 3$ and $n- 1 \ge \Delta \ge \delta \ge 2$, we note that
\[
\frac{2}{3}n - \frac{2}{3}(\Delta - 2) < \left( \frac{2\delta-1}{3\delta-2} \right) n - \frac{\Delta - 2}{3\delta-2}. \1
\]

Hence, we have shown that Dominator's greedy strategy is guaranteed to complete the total isolation game in at most~$\left( \frac{2\delta-1}{3\delta-2} \right) n - \frac{\Delta - 2}{3\delta-2}$ moves, yielding the desired upper bound in the statement of the theorem. This completes the proof of Theorem~\ref{thm:delta-bound}.~\QED

\medskip
If $\delta = \Delta = 2$, then by Theorem~\ref{thm:delta-bound} we have $\igt(G) \le \left( \frac{2\delta-1}{3\delta-2} \right) n$.
If $\delta \ge 2$ and $\Delta \ge 3$, then by Theorem~\ref{thm:delta-bound} we have $\igt(G) < \left( \frac{2\delta-1}{3\delta-2} \right) n $.
Hence as an immediate consequence of Theorem~\ref{thm:delta-bound}, we obtain the following simplified upper bound on the game total isolation number.

\begin{corollary}
\label{cor:delta-bound-1}
If $G$ is a connected graph of order~$n$ with minimum degree~$\delta$ where $\delta \ge 2$, then
\[
\igt(G) \le \left( \frac{2\delta-1}{3\delta-2} \right) n,
\]
with strict inequality if the maximum degree $\Delta \ge 3$.
\end{corollary}

An analogous proof as in Theorem~\ref{thm:delta-bound} yields the following upper bound on the Staller-start total isolation game. For completeness we provide the proof details.

\begin{theorem}
\label{thm:delta-bound-Staller}
If $G$ is a connected graph of order~$n \ge 3$ with minimum degree~$\delta$ where $\delta \ge 2$, then
\[
\igtS(G) \le \left( \frac{2\delta-1}{3\delta-2} \right) n - \frac{(\delta - 1)(2\delta - 3)}{3\delta-2}.
\]
\end{theorem}
\proof
Let $G$ be a connected graph of order~$n \ge 3$ with minimum degree~$\delta$ where $\delta \ge 2$ and maximum degree~$\Delta$. Suppose that $\Delta = n-1$. Thus, $G$ contains a dominating vertex, that is, a vertex adjacent to every other vertex. In this case, if \St does not play a dominating vertex on her first move, then \D will do so, and complete the game in two moves. Hence, $\igtS(G) = 2$. Moreover since $\delta \le n - 1$, we have $\frac{2}{3}(n - \delta + 2) \ge 2$, and so $\igtS(G) \le \frac{2}{3}(n - \delta + 2)$. Hence we may assume that $\Delta \le n-2$. In particular, $n \ge 4$. If $n = 4$, then $G = C_4$ and $\igtS(G) = 2 < \frac{2}{3}(n - \delta + 2)$. Hence, we may further assume that $n \ge 5$.

As the game is played, let $m_{2i}$ be the number of new vertices marked by \D on his $i$th move, and let $m_{2i-1}$ be the number of new vertices marked by \St on her $i$th move. Upon completion of the game, all vertices are marked implying that
\begin{equation}
\label{Eq1:St}
n = |M| = \sum_{i=1}^{\lfloor \frac{t}{2} \rfloor} m_{2i} + \sum_{i=1}^{\lceil \frac{t}{2} \rceil} m_{2i-1}.
\end{equation}

In what follows, let \D play according to a greedy strategy that on each of his moves he plays a vertex that results in a maximum increase in the number of marked vertices. Before Staller's first move, all vertices are unmarked, and so $U = V(G)$ and $S = M = \emptyset$. Let $t$ denote the number of moves played upon completion of the game when \D adopts his greedy strategy.

We consider Stage~$1$ and Stage~$2$ in the game as defined in the proof of Theorem~\ref{thm:delta-bound}. Suppose that Stage~1 consists of $k_1$ moves of Dominator. Since $n \ge 5$ and $\Delta \le n - 2$ by assumption, we note that $k_1 \ge 1$. Each of Dominator's $k_1$ moves in Stage~1 mark at least two new vertices, and so $m_{2i} \ge 2$ for all $i \in [k_1]$. Staller's first move marks at least~$\delta$ vertices and each of her remaining moves mark at least one new vertex.

\begin{claim}
\label{c:claim-3}
If the game is complete in Stage~$1$, then $t < \frac{2}{3}(n - \delta + 2)$.
\end{claim}
\proof
If the game is complete immediately after Dominator plays his $k_1^{\rm th}$ move, then \St plays $k_1$ moves in total, and so $t = 2k_1$ and, by Equation~(\ref{Eq1:St}), we have
\[
n = \sum_{i=1}^{k_1} m_{2i} + \sum_{i=1}^{k_1} m_{2i-1} \ge 2k_1 + (\delta + k_1 - 1) = 3k_1 + \delta - 1,
\]
and so
\[
t = 2k_1 = \frac{2}{3}(n - \delta + 1) < \frac{2}{3}(n - \delta + 2), \1
\]
and so the desired upper bound in the statement of the claim holds. Hence, we may assume that $t \ge 2k_1 + 1$. Since the game is complete in Stage~$1$, we infer that \St plays $k_1 + 1$ moves and $t = 2k_1 + 1$. By Equation~(\ref{Eq1:St}),  we have $n \ge 2k_1 + (\delta + k_1) = 3k_1 + \delta$, and in this case
\[
t = 2k_1 + 1 \le \frac{2}{3}(n  - \delta) + 1  < \frac{2}{3}(n - \delta + 2), \1
\]
and so, once again, the desired upper bound in the statement of the claim holds.~\smallqed

\medskip
Suppose, next, that the game is not complete in Stage~$1$, and so the game enters Stage~$2$.  Let Stage~$1$ consists of $k_1$ moves of Dominator. In Stage~$2$ of the game, every playable move by either player marks exactly one new vertex. Suppose that $k_2$ vertices in total are played in Stage~$2$ of the game, and so
\begin{equation}
\label{Eq3:St}
t = 2k_1 + k_2 + 1.
\end{equation}

Claim~\ref{c:claim-2} in the proof of Theorem~\ref{thm:delta-bound} holds exactly as before. Defining the sets $A$ and $B$ as in the proof of Theorem~\ref{thm:delta-bound}, we have
\begin{equation}
\label{Eq4:St}
|B| \ge \delta k_2 \hspace*{0.5cm} \mbox{and} \hspace*{0.5cm} |A| \ge 2k_1 + 1 - k_2.
\end{equation}

By Inequality~(\ref{Eq4:St}), and recalling that $|U| = k_2$, we have
\[
\begin{array}{lcl}
k_2 & = & n - |A| - |B| \\ \1
& \le & n - (2k_1 + 1 - k_2) - \delta  k_2 \\ \1
& = & n - (2k_1 + 1) - (\delta - 1)k_2,
\end{array}
\]
and so $(\delta - 1)k_2 \le n - (2k_1 + 1 + k_2) = n - t$, or, equivalently,
\begin{equation}
\label{Eq6:St}
k_2 \le \frac{n-t}{\delta - 1}.
\end{equation}

Immediately upon completion of Stage~$1$ of the game, we have
\[
|M| = \sum_{i=1}^{k_1} m_{2i} + \sum_{i=1}^{k_1+1} m_{2i-1} \ge 2k_1 + (\delta + k_1) = 3k_1 + \delta.
\]

Recall that $t = 2k_1 + 1 + k_2$ by Equation~(\ref{Eq3:St}), and so $k_1 = \frac{1}{2}(t - k_2 - 1)$. From these observations and by Inequality~(\ref{Eq6:St}), we have
\[
\begin{array}{lcl}
n & = & |M| + |U| \\ \1
& \ge & (3k_1 + \delta) + k_2 \\ \1
& = & \frac{3}{2}(t - k_2 - 1) + k_2 + \delta \\ \1
& = & \frac{3}{2}t - \frac{1}{2}k_2 + \delta - \frac{3}{2} \\ \1
& \ge & \frac{3}{2}t - \frac{n-t}{2(\delta - 1)} + \delta - \frac{3}{2},
\end{array}
\]
or, equivalently,
\begin{equation}
\label{Eq8:St}
\left( \frac{2\delta - 1}{2\delta - 2} \right) n \ge \left( \frac{3\delta - 2}{2\delta - 2} \right) t + \delta - \frac{3}{2}.
\end{equation}

Solving for $t$ yields the inequality
\begin{equation}
\label{Eq9:St}
t \le \left( \frac{2\delta-1}{3\delta-2} \right) n - \frac{(\delta - 1)(2\delta - 3)}{3\delta-2}. \1
\end{equation}

Recall that $t$ denotes the number of moves played in the game when \D adopts his greedy strategy. By Claim~\ref{c:claim-3}, if the game is complete in Stage~$1$, then $t < \frac{2}{3}(n - \delta + 2)$. If the game is not complete in Stage~$1$ and therefore enters Stage~$2$, then by Inequality~(\ref{Eq9:St}) we have $t \le \left( \frac{2\delta-1}{3\delta-2} \right) n - \frac{(\delta - 1)(2\delta - 3)}{3\delta-2}$. For $n \ge 3$ and $n- 1 \ge\delta \ge 2$, we note that
\[
\frac{2}{3}(n - \delta + 2) \le  \left( \frac{2\delta-1}{3\delta-2} \right) n - \frac{(\delta - 1)(2\delta - 3)}{3\delta-2}. \1
\]

Hence, we have shown that Dominator's greedy strategy is guaranteed to complete the total isolation game in at most~$\left( \frac{2\delta-1}{3\delta-2} \right) n - \frac{(\delta - 1)(2\delta - 3)}{3\delta-2}$ moves, yielding the desired upper bound in the statement of the theorem. This completes the proof of Theorem~\ref{thm:delta-bound-Staller}.~\QED

\medskip
If $\delta = 2$, then by Theorem~\ref{thm:delta-bound-Staller} we have $\igtS(G) \le \left( \frac{2\delta-1}{3\delta-2} \right) n - \frac{1}{4}$.

If $\delta \ge 3$, then by Theorem~\ref{thm:delta-bound-Staller} we have $\igtS(G) < \left( \frac{2\delta-1}{3\delta-2} \right) n - \frac{1}{4}$.

Hence as an immediate consequence of Theorem~\ref{thm:delta-bound-Staller}, we obtain the following simplified upper bound on the Staller-start game total isolation number.

\begin{corollary}
\label{cor:delta-bound-Staller-1}
If $G$ is a connected graph of order~$n \ge 3$ with minimum degree~$\delta$ where $\delta \ge 2$, then
\[
\igtS(G) \le \left( \frac{2\delta-1}{3\delta-2} \right) n - \frac{1}{4},
\]
with strict inequality if $\delta \ge 3$.
\end{corollary}

For $\delta \ge 2$, we note that $\frac{2\delta-1}{3\delta-2} \le \frac{3}{4}$ (with strict inequality if $\delta \ge 3$). Hence, Corollary~\ref{cor:delta-bound-1} and Corollary~\ref{cor:delta-bound-Staller-1} yield the following upper bounds on the (Dominator-start) game total isolation number and Staller-start game total isolation number of a graph with minimum degree at least~$2$.

\begin{corollary}
\label{cor:delta-bound-2}
If $G$ is a connected graph of order~$n \ge 3$ with $\delta(G) \ge 2$, then the following hold.  \\ [-24pt]
\begin{enumerate}
\item[{\rm (a)}] $\igt(G) \le \frac{3}{4}n$. \2
\item[{\rm (b)}] $\igtS(G) \le \frac{3}{4}n - \frac{1}{4}$.
\end{enumerate}
\end{corollary}

\subsection{Graphs with diameter two}

Let $G$ be a connected graph of order~$n \ge 3$ with minimum degree~$\delta$ where $\delta \ge 2$. By Corollary~\ref{cor:delta-bound-1}, $\igt(G) \le \left( \frac{2\delta-1}{3\delta-2} \right) n$ and by Corollary~\ref{cor:delta-bound-2}, $\igtS(G) \le \left( \frac{2\delta-1}{3\delta-2} \right) n - \frac{1}{4}$. Hence if $\delta$ tends to infinity, then the upper bound on $\igt(G)$ tends to $\frac{2}{3}n$ and the upper bound on $\igtS(G)$ tends to $\frac{2}{3}n - \frac{1}{4}$.

As a consequence of our proof of Theorem~\ref{thm:delta-bound}, we show that the (Dominator-start) game total isolation number and Staller-start game total isolation numbers for connected graphs of order~$n \ge 3$ and diameter at most~$2$ is at most~$\frac{2}{3}n$. We remark that it is well-known that almost every graph has diameter~$2$ by a property of the Erd\H{o}s-R\'{e}nyi random graph model. Hence by Theorem~\ref{thm:diam2} below, almost every graph $G$ satisfies $\igt(G) \le \frac{2}{3}n$ and $\igtS(G) \le \frac{2}{3}n$

\begin{theorem}
\label{thm:diam2}
If $G$ is a connected graph of order~$n \ge 3$ with $\diam(G) \le 2$, then
\[
\igt(G) \le \frac{2}{3}n \hspace*{0.5cm} \mbox{and} \hspace*{0.5cm} \igtS(G) \le \frac{2}{3}n.
\]
\end{theorem}
\proof
Let $G$ be a connected graph of order~$n\ge 3$ with $\diam(G) \le 2$. Further, let $G$ have minimum degree~$\delta$ where $\delta \ge 1$ and maximum degree~$\Delta$. If $G=K_n$, then $2 = \igt(G) \le \frac{2}{3}n$  and $2 = \igtS(G) \le \frac{2}{3}n$. Thus, we assume that $\diam(G) = 2$, which implies that $\Delta \ge 2$.

Suppose that $\delta = 1$. Let $v_1$ be a vertex of degree~$1$ in $G$ and let $v$ be its (unique) neighbor. Since $\diam(G) = 2$, we note that $V(G) = N_G[v]$, and so $\deg_G(v) = \Delta = n-1$. Thus, $G$ contains a dominating vertex, namely~$v$.
In the D-game, if \D plays the vertex~$v$ on his first move, he marks all vertices of $G$ different from~$v$. Since one additional move completes the game, we have $\igt(G)  = 2 \le \frac{2}{3}n$.
In the S-game, if \St plays the dominating vertex~$v$, then one additional move of \D completes the game, and if \St does not play the vertex~$v$ on her first move, then \D plays the vertex~$v$ on his first move, and once again the game is complete in two moves. Hence, $\igtS(G) = 2 \le \frac{2}{3}n$.
Hence we may assume that $\delta \ge 2$.

Suppose that $\Delta = 2$. Since $\diam(G) = 2$, we infer that $G \in \{C_4,C_5\}$. Since $\igt(C_4) = 2$ and $\igt(C_5) = 3$, we have that $\igt(G) < \frac{2}{3}n$. Moreover since $\igtS(C_4) = 2$ and $\igtS(C_5) = 2$, we have that $\igt(G) < \frac{2}{3}n$. Hence we may further assume that $\Delta \ge 3$.

Suppose firstly that we are in the D-game. In this case, we adopt the notation and proof of Theorem~\ref{thm:delta-bound}. In particular, $t$ denotes the number of moves played in the D-game when \D adopts his greedy strategy, and so $\igt(G) \le t$. We show that $t \le \frac{2}{3}n$. By our earlier assumptions, $\delta \ge 2$ and $\Delta \ge 3$. By Claim~\ref{c:claim-2} in the proof of Theorem~\ref{thm:delta-bound} we have that $U \subseteq S$ and $U$ is a packing in $G$. If $|U| \ge 2$, then $\diam(G) \ge 3$, a contradiction. Hence, $k_2 = |U| \le 1$.
Suppose that $k_2 = 0$. Thus, the game is completed in Stage~$1$, and $t = 2k_1 - 1$ or $t = 2k_1$. If $t = 2k_1 - 1$, then as shown in the proof of Stage~$1$ of Theorem~\ref{thm:delta-bound} we have $t < \frac{2}{3}n$. Hence, we may assume that $t = 2k_1$. Recall that $m_{2i-1} \ge 2$ for all $i \in [k_1]$ and $m_{2i} \ge 1$ for all $i \in [k_1]$. As shown in the proof of Theorem~\ref{thm:delta-bound} following Dominator's greedy strategy we have $m_1 \ge \Delta$, and so $m_1 \ge 3$. Thus,
\[
n = \left( m_1 + \sum_{i=2}^{k_1} m_{2i-1} \right)  + \left( \sum_{i=1}^{k_1} m_{2i} \right)
\ge (2k_1 + 1) + k_1 = 3k_1 + 1,
\]
and so, $t = 2k_1 \le \frac{2}{3}(n-1) < \frac{2}{3}n$. Hence we may assume that $k_2 = 1$, for otherwise the desired upper bound follows. In this case, we have $t = 2k_1 + k_2 = 2k_1 + 1$, and so
\[
n = \left( m_1 + \sum_{i=2}^{k_1} m_{2i-1} \right)  + \left( \sum_{i=1}^{k_1} m_{2i} \right) + k_2
\ge (2k_1 + 1) + k_1 + 1 = 3k_1 + 2,
\]
implying that $t = 2k_1 + 1 \le \frac{2}{3}(n-2) + 1 < \frac{2}{3}n$.

Suppose secondly that we are in the S-game. In this case, we adopt the notation and proof of Theorem~\ref{thm:delta-bound-Staller}. In particular, $t$ denotes the number of moves played in the S-game when \D adopts his greedy strategy, and so $\igt(G) \le t$. By our earlier assumptions, $\delta \ge 2$ and $\Delta \ge 3$.  We show that $t \le \frac{2}{3}n$. As in the case of the D-game we have that $U \subseteq S$. Further, the set $U$ is a packing in $G$ and $k_2 = |U| \le 1$.
Suppose that $k_2 = 0$. Thus, the game is completed in Stage~$1$, and $t = 2k_1$ or $t = 2k_1 + 1$. If $t = 2k_1$, then as shown in the proof of Theorem~\ref{thm:delta-bound-Staller} we have $t \le \frac{2}{3}(n - \delta + 1) < \frac{2}{3}n$ noting that $\delta \ge 2$. Hence, we may assume that $t = 2k_1 + 1$. Recall that $m_1 \ge \delta \ge 2$ and $m_{2i-1} \ge 1$ for all $i \in [k_1 + 1] \setminus \{1\}$, and $m_{2i} \ge 2$ for all $i \in [k_1]$. Thus,
\[
n = \left( \sum_{i=1}^{k_1} m_{2i} \right)  + \left( m_1 + \sum_{i=2}^{k_1+1} m_{2i-1} \right)
\ge 2k_1 + (2 + k_1) = 3k_1 + 2,
\]
and so, $t = 2k_1 + 1 \le \frac{2}{3}(n-2) + 1 < \frac{2}{3}n$. Hence we may assume that $k_2 = 1$, for otherwise the desired upper bound follows. In this case, we have $t = 2k_1 + 1 + k_2 = 2k_1 + 2$, and so
\[
n = \left( \sum_{i=1}^{k_1} m_{2i} \right)  + \left( m_1 + \sum_{i=1}^{k_1 + 1} m_{2i-1} \right) + k_2
\ge 2k_1 + (2 + k_1) + 1 = 3k_1 + 3,
\]
implying that $t = 2k_1 + 2 \le \frac{2}{3}(n-3) + 2 = \frac{2}{3}n$.~\QED

\section{General upper bounds for connected graphs}
\label{S:min-degree-1}

We consider next general upper bounds on the (Dominator-start) game total isolation number and Staller-start game total isolation number of a connected graph.

\begin{theorem}
\label{thm:min-deg-1}
If $G$ is a connected graph of order~$n \ge 3$, then $\igt(G) < \frac{5}{6}n$.
\end{theorem}
\proof
Let $G$ be a connected graph of order~$n \ge 3$ with minimum degree~$\delta$ where $\delta \ge 1$ and maximum degree~$\Delta$. If $\delta \ge 2$, then the desired bound follows from Corollary~\ref{cor:delta-bound-1}. If $\diam(G) \le 2$, then the desired bound follows from Theorem~\ref{thm:diam2}. Hence we may assume that $\delta = 1$ and that the connected graph $G$ satisfies $\diam(G) \ge 3$. In what follows we adopt the notation and proof of Theorem~\ref{thm:delta-bound}. In particular, $t$ denotes the number of moves played in the game when \D adopts his greedy strategy, and so $\igt(G) \le t$.

We now modify Dominator's strategy as follows. As in the proof of Theorem~\ref{thm:delta-bound},  \D plays according to a greedy strategy, and so on each of his moves, \D plays a vertex that marks as many previously unmarked vertices as possible. Subject to \D playing according to his greedy strategy, \D selects a non-leaf on each of his moves whenever possible. Thus if \D has two or more greedy moves available to play (each marking the same number of previously unmarked vertices), then \D selects among all such greedy moves a non-leaf if at all possible. We show that with the modified strategy, we have $t \le \frac{5}{6}n$.

Consider the sets $S$, $M$ and $U$ of selected, marked and unmarked vertices immediately after Stage~1 and before the game enters Stage~2. By Claim~\ref{c:claim-2} in the proof of Theorem~\ref{thm:delta-bound} we have that $U \subseteq S$ and $U$ is a packing in $G$. Recall that $k_2 = |U|$. If $k_2 = 0$, then the game is completed in Stage~$1$, and $t = 2k_1 - 1$ or $t = 2k_1$. In both cases, as shown in the proof of Stage~$1$ of Theorem~\ref{thm:delta-bound} we have $t \le \frac{2}{3}n$. Hence, we may assume that $k_2 \ge 1$.

Adopting our earlier notation in the proof of Theorem~\ref{thm:delta-bound}, let $U = \{u_1,\ldots,u_{k_2}\}$ and let $B_i = N_G(u_i)$ for $i \in [k_2]$. As before, let $B = N_G(U)$.  Renaming vertices in $U$ if necessary, we may assume $\deg_G(u_1) \le \deg_G(u_2) \le \cdots \le \deg_G(u_{k_2})$. Let $U_1$ be the subset of vertices in $U$ of degree~$1$ in $G$, and let $U_2 = U \setminus U_1$, and so every vertex in $U_2$ has degree at least~$2$ in $G$. Further, let $|U_i| = n_i$ for $i \in [2]$. If $n_1 = 0$, then $U = U_2$ and
\begin{equation}
\label{Eq:min-deg1-Eq1}
|B| = \sum_{i=1}^{k_2} \deg_G(u_i) \ge 2k_2.
\end{equation}

Recall that by Inequality~(\ref{Eq5}) in the proof of Theorem~\ref{thm:delta-bound} we have
$|A| \ge 2k_1 - k_2$. Together with Inequality~(\ref{Eq:min-deg1-Eq1}) this yields
\[
\begin{array}{lcl}
k_2 \, = \, |U| & = & n - |A| - |B| \\ \1
& \le & n - (2k_1 - k_2) - 2k_2 \\ \1
& = & n - (2k_1 + k_2) \\ \1
& = & n - t.
\end{array}
\]
Thus,
\[
\begin{array}{lcl}
n \, = \, |M| + |U|  & \ge & 3k_1 + k_2 \\ \1
& = & \frac{3}{2}(t - k_2) + k_2 \\ \1
& = & \frac{3}{2}t - \frac{1}{2}k_2 \\ \1
& \ge & \frac{3}{2}t - \frac{1}{2}(n-t),
\end{array}
\]
or, equivalently, $t \le \frac{3}{4}n$. Hence we may assume that $n_1 \ge 1$.

\begin{unnumbered}{Claim~A}
$k_1 \ge n_1$.
\end{unnumbered}
\proof
Let $u$ be an arbitrary vertex in $U_1$. Thus, $\deg_G(u) = 1$. Let $v$ be the (unique) neighbor of~$u$. By assumption, $\diam(G) \ge 3$, implying that $v$ has a neighbor, say~$w$, of degree at least~$2$ in $G$. By our earlier observations, $U_1 \subseteq U \subseteq S$ and $U$ is a packing in $G$. Moreover, $k_2 = |U|$. Recall that $S$, $M$ and $U$ denote the sets of selected, marked and unmarked vertices immediately after Stage~1 and before the game enters Stage~2. In particular, we note that all vertices in $U$ were played in Stage~$1$ of the game.

We show that the vertex~$u \in U_1$ was played by Staller (in Stage~$1$ of the game). Suppose, to the contrary, that \D played the vertex~$u$ in Stage~$1$ of the game. Let $S'$, $M'$, and $U'$ denote the sets of selected, marked and unmarked vertices immediately before vertex~$u$ is played. Then, $u \notin S'$ and $v \in U'$.  Also, since $u$ is unmarked at the end of Stage~$1$, it follows that $u \in U'$. However if \D had played the vertex~$w$ instead of the vertex~$u$, then he would have marked all the new vertices that were marked when playing the vertex~$u$, except possibly for the vertex~$w$ (which may become marked when the vertex~$u$ is played), and in addition he would have marked the vertex~$u$. Hence, \D would have marked at least as many vertices by playing the vertex~$w$ instead of playing the vertex~$u$. However since $u$ is a leaf and $w$ is a non-leaf, according to Dominator's modified strategy he would have played the vertex~$w$, a contradiction. Hence the vertex~$u$ must have been played by Staller in Stage~$1$ of the game. Moreover when the vertex~$u$ was played by Staller, Dominator's previous (greedy) move necessary marked at least two new vertices, noting that playing the vertex~$w$ would mark at least two new vertices, namely~$u$ and~$v$.

The vertices in the set $U_1$ were therefore all selected by Staller. Since every move of Staller that played a vertex in $U_1$ was played in Stage~1 of the game and since \St played a total of $k_1$ moves in Stage~1 of the game, we infer that $n_1 \le k_1$.~\smallqed

\begin{unnumbered}{Claim~B}
$k_2 \le \frac{1}{2}(n - k_1)$.
\end{unnumbered}
\proof
By assumption, $n_1 \ge 1$. Thus,
\begin{equation}
\label{Eq:min-deg1-Eq2}
|B| = \sum_{i=1}^{n_1} |B_i| + \sum_{i=n_1 + 1}^{k_2} |B_i|  \ge n_1 + 2(k_2 - n_1) = 2k_2 - n_1.
\end{equation}

Recall that by Inequality~(\ref{Eq5}) in the proof of Theorem~\ref{thm:delta-bound} we have
$|A| \ge 2k_1 - k_2$. By Claim~A, we have $k_1 \ge n_1$. Together with Inequality~(\ref{Eq:min-deg1-Eq2}) these observations yield the inequality chain given by
\[
\begin{array}{lcl}
k_2 \, = \, |U| & = & n - |A| - |B| \1 \\
& \le & n - (2k_1 - k_2) - (2k_2 - n_1) \1 \\
& = & n - 2k_1 - k_2 + n_1 \1 \\
& \le & n - k_1 - k_2,
\end{array}
\]
or, equivalently, $k_2 \le \frac{1}{2}(n - k_1)$.~\smallqed

\medskip
Recall that $n \ge \frac{3}{2}t - \frac{1}{2}k_2$. Hence by Claim~B, we infer that $n \ge \frac{3}{2}t - \frac{1}{4}(n - k_1)$, or, equivalently, $t \le \frac{5}{6}n - \frac{1}{6}k_1$. Since $k_1 \ge 1$, we have $t < \frac{5}{6}n$, as claimed. This completes the proof of Theorem~\ref{thm:min-deg-1}.~\QED

\begin{theorem}
\label{thm:min-deg-Staller-1}
If $G$ is a connected graph of order~$n \ge 3$, then $\igtS(G) \le \frac{5}{6}n$.
\end{theorem}
\proof
Let $G$ be a connected graph of order~$n \ge 3$ with minimum degree~$\delta$ where $\delta \ge 1$ and maximum degree~$\Delta$. If $\delta \ge 2$, then by Corollary~\ref{cor:delta-bound-2} we have $\igtS(G) \le \frac{3}{4}n - \frac{1}{4} < \frac{5}{6}n$.
If $\diam(G) \le 2$, then the desired bound follows from Theorem~\ref{thm:diam2}. Hence we may assume that $\delta = 1$ and that the connected graph $G$ satisfies $\diam(G) \ge 3$. In what follows we adopt the notation and proof of Theorem~\ref{thm:delta-bound-Staller}. In particular, $t$ denotes the number of moves played in the S-game when \D adopts his greedy strategy, and so $\igt(G) \le t$.

We now modify Dominator's strategy as follows. As in the proof of Theorem~\ref{thm:diam2},  \D plays according to a greedy strategy, and so on each of his moves, \D plays a vertex that marks as many previously unmarked vertices as possible. Subject to \D playing according to his greedy strategy, \D selects a non-leaf on each of his moves whenever possible. Thus if \D has two or more greedy moves available to play (each marking the same number of previously unmarked vertices), then \D selects among all such greedy moves a non-leaf if at all possible. We show that with the modified strategy, we have $t \le \frac{5}{6}n$.

Consider the sets $S$, $M$ and $U$ of selected, marked and unmarked vertices immediately after Stage~1 and before the game enters Stage~2. By Claim~\ref{c:claim-2} in the proof of Theorem~\ref{thm:delta-bound} we have that $U \subseteq S$ and $U$ is a packing in $G$. Recall that $k_2 = |U|$. If $k_2 = 0$, then the game is completed in Stage~$1$, and $t = 2k_1$ or $t = 2k_1 + 1$. In both cases, as shown in the proof of Stage~$1$ of Theorem~\ref{thm:delta-bound-Staller} we have $t \le \frac{2}{3}n + \frac{1}{3}$ noting that $\delta = 1$. Hence, we may assume that $k_2 \ge 1$.

We now adopt our earlier notation for the set $A$ and $B$ as in the proof of Theorem~\ref{thm:delta-bound-Staller} and for the sets $U$, $U_1$ and $U_2$  in the proof of Theorem~\ref{thm:min-deg-1}. In particular, $|U_i| = n_i$ for $i \in [2]$. If $n_1 = 0$, then $U = U_2$ and, as in the proof of Theorem~\ref{thm:delta-bound-Staller}, we infer in this case that
\begin{equation}
\label{Eq14:St}
|B| \ge 2 k_2 \hspace*{0.5cm} \mbox{and} \hspace*{0.5cm} |A| \ge 2k_1 + 1 - k_2.
\end{equation}
By Inequality~(\ref{Eq14:St}), and recalling that $|U| = k_2$, we have
\[
\begin{array}{lcl}
k_2 \, = \, |U| & = & n - |A| - |B| \\ \1
& \le & n - (2k_1 + 1 - k_2) - 2k_2 \\ \1
& = & n - (2k_1 + 1 + k_2) \\ \1
& = & n - t.
\end{array}
\]
Thus,
\[
\begin{array}{lcl}
n \, = \, |M| + |U|  & \ge & (3k_1 + 1) + k_2 \\ \1
& = & \frac{3}{2}(t - k_2 - 1) + k_2 + 1 \\ \1
& = & \frac{3}{2}t - \frac{1}{2}k_2 - \frac{1}{2}\\ \1
& \ge & \frac{3}{2}t - \frac{1}{2}(n-t) - \frac{1}{2} ,
\end{array}
\]
or, equivalently, $t \le \frac{3}{4}n + \frac{1}{4}$. Hence we may assume that $n_1 \ge 1$. Proceeding now exactly as in the proof of Theorem~\ref{thm:min-deg-1} (see Claim~A), we have $k_1 \ge n_1$ and $|B| \ge 2k_1 - n_1$. We therefore infer that
\[
\begin{array}{lcl}
k_2 \, = \, |U| & = & n - |A| - |B| \\ \1
& \le & n - (2k_1 + 1 - k_2) - (2k_2 - n_1) \\ \1
& = & n - k_1 - k_2 - 1 + (n_1 - k_1) \\ \1
& \le & n - k_1 - k_2 - 1,
\end{array}
\]
or, equivalently, $k_2 \le \frac{1}{2}(n - k_1 - 1)$. As shown earlier, we have $n \ge \frac{3}{2}t - \frac{1}{2}k_2 - \frac{1}{2}$. Hence, $n \ge \frac{3}{2}t - \frac{1}{2}k_2 - \frac{1}{2} \ge \frac{3}{2}t - \frac{1}{4}(n - k_1 - 1) - \frac{1}{2}$, or, equivalently, $t \le \frac{5}{6}n - \frac{1}{6}(k_1 + 1) + \frac{1}{3}$. Since $k_1 \ge 1$, we have $t \le \frac{5}{6}n$, as claimed. This completes the proof of Theorem~\ref{thm:min-deg-Staller-1}.~\QED

\section{Closing comments}

We pose the following $\frac{2}{3}$-conjecture for the game total isolation number.

\begin{conjecture}
\label{conj:1}
If $G$ is a graph of order~$n$ and every component of $G$ has order at least~$3$, then $\igt(G) \le \frac{2}{3}n$.
\end{conjecture}

We remark that if Conjecture~\ref{conj:1} is true, then the $\frac{2}{3}n$-upper bound on the game total isolation number is best possible. As an example, let $G$ be a graph of order~$n$ where every component is a path  $P_3$- or $P_6$-component or a cycle $C_3$- or $C_6$-component. Staller's strategy is to reply to every move of \D by playing in the component that \D played in on his previous move. When \D plays in a $P_3$- or $C_3$-component,  \St responds with a playable vertex in that component that marks all vertices. When \D plays in a $P_6$- or $C_6$-component, \St responds by playing a playable vertex in that component following her optimal strategy, thereby ensuring that four moves (two by \D and two by Staller) are played in such components. (Note that if \D plays a vertex~$v$ from a $C_6$-component $C$, then \St plays the vertex at distance~$3$ from~$v$ on the cycle $C$, thereby forcing two additional moves played on~$C$.) Upon completion of the game, two moves are played in every $P_3$- and $C_3$-component and four moves are played in every $P_6$- and $C_6$-component, implying that $\igt(G) \ge \frac{2}{3}n$. However if Conjecture~\ref{conj:1} is true, then $\igt(G) \le \frac{2}{3}n$. Consequently, $\igt(G) = \frac{2}{3}n$.


\medskip
\begin{thebibliography}{99}



\bibitem{Bo-25}
P.~Borg, Proof of a conjecture on isolation of graphs dominated by a vertex. \textit{Discrete Appl.\ Math.} \textbf{371} (2025), 247--253.

\bibitem{BoCh-18}
R. Boutrig and M. Chellali, Total vertex-edge domination. \textit{Int. J. Comput. Math.} \textbf{95} (2018), 1820--1828.

\bibitem{BoGo-24}
G.~Boyer and W.~Goddard, Disjoint isolating sets and graphs with maximum isolation number. \textit{Discrete Appl.\ Math.} \textbf{356} (2024), 110--116.

\bibitem{BoGoHe-25a}
G.~Boyer, W.~Goddard, and M. A.~Henning, On total isolation in graphs. \textit{Aequationes Math.} \textbf{99} (2025), 623--633.

\bibitem{BoGoHe-25b}
G.~Boyer, W.~Goddard, and M. A.~Henning, Bounds on the total isolation number of (regular) bipartite graphs, manuscript, 2025.

\bibitem{BrDrJoKuRa-24}
B.\ Bre\v sar, T.\ Dravec, D. P.\ Johnston, K.\ Kuenzel, and D. F.\ Rall, Isolation game on graphs.
\url{arXiv:2409.14180} (2024).


\bibitem{BrHeKlRa-21}
B. Bre\v{s}ar, M. A. Henning, S. Klav\v zar, and D. F. Rall, \emph{Domination Games Played on Graphs.} \emph{SpringerBriefs in Mathematics}. \emph{Springer, Cham}, (2021) x + 122 pp.


\bibitem{BrKlRa-10}
B.~Bre{\v{s}}ar, S.~Klav{\v{z}}ar, and D.~F.~Rall, Domination game and an imagination strategy. \textit{SIAM J. Discrete Math.} \textbf{24} (2010), 979--991.


\bibitem{BuDrHeKl-17}
Cs. Bujt\'as, T. Dravec, M. A. Henning, and S. Klav\v zar, Bounds on the game isolation number and exact values for paths and cycles. \url{arXiv:2507.08503} (2025).


\bibitem{CaAnWu-25}
Y. Cao, X. An, Xinhui, and B. Wu, Total isolation of $k$-cliques in a graph. \textit{Discrete Math.} \textbf{348} (2025), Paper No. 114689.

\bibitem{CaHa-17} Y.~Caro and A.~Hansberg, Partial domination---the isolation number of a graph. \textit{Filomat} \textbf{31} (2017), 3925--3944.

\bibitem{DaLeSoVa-21} A. Dapena, M. Lema\'{n}ska, M. J. Souto-Salorio, and F. J. Vazquez-Araujo, Isolation number versus domination number of trees. \textit{Mathematics} \textbf{9} (2021), Paper 1325.

\bibitem{GoHe-25}
W. Goddard and M. A. Henning, On the isolation number of graphs with minimum degree four.  \url{arXiv:2508.21551} (2025).


\bibitem{HaHeHe-23}
T. W. Haynes, S. T. Hedetniemi, and M. A. Henning, \emph{Domination in Graphs: Core Concepts} Series: Springer Monographs in Mathematics, Springer, Cham, 2023. xx + 644 pp.


\bibitem{HeYe-book}
M. A. Henning and A. Yeo, \emph{Total domination in graphs.} Series: Springer Monographs in Mathematics, Springer, Cham, New York, 2013. xiv + 178 pp.

\bibitem{KiWeZa-13}
W.~B.~Kinnersley, D.~B.~West, and R.~Zamani, Extremal problems for game domination number. \textit{SIAM J. Discrete Math.} \textbf{27} (2013), 2090--2107.

\bibitem{LeMoSo-24}
M.~Lema\'{n}ska, M.~Mora, and M. J.~Souto-Salorio, Graphs with isolation number equal to one third of the order. \textit{Discrete Math.} \textbf{347} (2024), Paper 113903.

\bibitem{Pe-86} J. W. Peters, Theoretical and algorithmic results on domination and connectivity. Ph.D. Thesis, Clemson University, Clemson (1986).

\bibitem{Zy-19} P. \.{Z}yli\'{n}ski, Vertex-edge domination in graphs. Aequationes Math. 93 (2019), 735--742.


\end{thebibliography}
\end{document}